# Solvable cases of optimal control problems for integral equations

S. A. Belbas, Mathematics Department, University of Alabama, Tuscaloosa, AL 35487-0350, USA; e-mail sbelbas@gmail.com

W. H. Schmidt, Institut für Mathematik und Informatik, Universität Greifswald, Germany.

Abstract. We present a number of cases of optimal control of Volterra and Fredholm integral equations that are solvable in the sense that the problem can be reduced to a solvable integral equation. This is conceptually analogous to the role of the Riccati differential system in the optimal control of ordinary differential equations.

1.     Introduction.

In the ordinary theory of optimal control (i.e. optimal control of systems governed by ordinary differential equations), a prominent role is played by problems for which the optimality equations (Hamiltonian equations and Pontryagin's maximum and/or dynamic programming equations), be reduced to solvable equations. The archetypical example of this type of reduction is the Riccati differential system for optimal control of linear ODE control system with a quadratic performance functional. (A contemporary account of Riccati equations in optimal control is contained in [Z].) This is a solvable optimal control problem, in the sense that the Riccati differential system can (under suitable, but quite general assumptions) be solved by standard well-behaved and well-analyzed algorithms for systems of ordinary differential equations. Once the solution of the Riccati differential system has been calculated, the corresponding optimal control is determined in the form of a linear causal feedback law.
In the above discussion, we do not consider the method of dynamic programming, since that has no useful counterpart in the optimal control of integral equations.
In this paper, we are interested in discovering analogous cases of solvability for the optimal control of systems governed by integral equations. We remark that examples of applied problems that lead to optimal control questions for integral equations can be found, e.g., in [AHS], [BS]. For integral equations of either Fredholm or Volterra types, there is no possibility (in general) of finding an optimal control in the form of a causal feedback law. Even if the integral equation that describes the system dynamics is causal, i.e. a Volterra integral equation, an optimal control cannot, in general, be obtained in causal form. Indeed, for a controlled Volterra integral equation

$$y(t) = y_0(t) + \int_0^t f(t,s,y(s),u(s))\,ds \qquad \text{--- (1.1)}$$

with cost functional (which is to be minimized)

$$J := \int_0^T F(t,y(t),u(t))\,dt \qquad \text{--- (1.2)}$$

the Hamiltonian is

$$H(t,y,\psi(\cdot),u) = F(t,y,u) + \int_t^T \psi(s) f(s,t,y,u)\,ds \qquad \text{--- (1.3)}$$

where $\psi$ is the co-state and its values are co-vectors (row vectors, if the state $y$ is represented as a column vector). The co-state satisfies the integral equation

$$\psi(t) = \nabla_y H(t,y,\psi(\cdot),u) \qquad \text{--- (1.4)}$$

(here we follow the convention of representing the gradient as a co-vector). For details of the mathematical theory and proofs, see [S1, S2].

An extremum principle similar to Pontryagin's principle says that an optimal control $u^*(\cdot)$ minimizes, a.e. in $t$, $H(t, y(t), \psi(\cdot), u)$.
`

Eq. (1.4) is a Volterra integral equation in reverse time, equation over [t, T], and it is linear in $\psi$, thus the solution can (in general) be expressed via a resolvent kernel $R$ as

$$\psi(t) = \nabla_y F(t, y(t), u(t)) + \int_t^T \nabla_y F(s, y(s), u(s)) R(s,t) \, ds \quad . \tag{1.5}$$

The resolvent kernel $R(s, t)$ is (among other dependencies) a functional of the restrictions of $y(\cdot)$, $u(\cdot)$ to the interval $(t, T)$. We use the notation $[y]_{t,T}$, $[u]_{t,T}$ to denote these restrictions and, for every function $g$ with domain (0, T), we denote by $[g]_{a,b}$ the restriction of the function $g$ to the interval $(a, b)$. Then the resolvent kernel depends also on $[y]_{t,T}$, $[u]_{t,T}$:

$$R(s,t) = R(s,t; [y]_{t,T}, [u]_{t,T}) \quad . \tag{1.6}$$

Consequently $\psi(t)$ is a functional of $[y]_{t,T}, [u]_{t,T}$. It follows from (1.5) that also $[\psi]_{t,T}$ is a function of $[y]_{t,T}, [u]_{t,T}$. Thus a minimizer $u_*(t)$ (assuming its existence) will be a (in general, set-valued) function of $[y]_{t,T}$. Substituting one selection of this set-valued function into the state dynamics (1.1) shows that the right-hand side of the Volterra integral equation (1.1) will depend on both $[y]_{0,t}$ through the dependence of $f$ on $y$) and on $\{[y]_{s,T} : 0 \le s \le t\}$ (through the dependence of $u(s)$ on $[y]_{s,T}$), thus on $[y]_{0,T}$. Therefore, at the end, we have an integral – functional equation with right-hand side, for each $t$, depending on $[y]_{0,T}$. This shows the Fredholm nature of an optimal problem for state dynamics given by a Volterra integral equation.

For comparison with previous results, we mention that, in the particular case a Volterra integral equation arising from the integral formulation of a differential equation, as in [PY], causal synthesis of an optimal control is possible. In the case of general Volterra integral equations with memory effects, as in the present paper, optimal controls cannot be expressed in the form of causal feedback.

3...

2.  Optimization of certain types of quadratic functionals.

Our exposition of quadratic control of linear integral equations will be facilitated by examining first some general quadratic optimization problems over $L^2$ spaces, independently of any control interpretations. In this section, we gather a few results on quadratic functionals of the form

$$E := \tfrac{1}{2}\iint_{G\times G} w^T(x) K_2(x,y) w(y) \, dx \, dy + \int_G \{\tfrac{1}{2} w^T(x) K_1(x) w(x) + r_0^T(x) w(x)\} \, dx \; .$$

--- (2.1)

Without loss of generality, the matrix-valued functions $K_1$ and $K_2$ can be assumed to be symmetric functions, in the following sense:

$$K_1^T(x) = K_1(x), \quad K_2^T(x,y) = K_2(y,x) \; .$$

--- (2.2)

This is the case because $K_1$ and $K_2$ can be replaced, in the definition of the functional $E$, by their respective symmetrizations $\tilde{K}_1$, $\tilde{K}_2$, defined by

$$\tilde{K}_1(x) := \tfrac{1}{2}[K_1(x) + K_1^T(x)], \quad \tilde{K}_2(x,y) := \tfrac{1}{2}[K_2(x,y) + K_2^T(y,x)] \; ,$$

--- (2.3)

without affecting the values of $E$.

We extend $E$ to a quadratic functional $\tilde{E}$ on the space $(L^2(G \mapsto \mathbb{R}^n))^2$ as follows:

$$\tilde{E} \equiv \tilde{E}(w,v) := \tfrac{1}{4}\iint_{G\times G} [w^T(x) \; v^T(y)] \tilde{K}(x,y) \begin{bmatrix} w(x) \\ v(y) \end{bmatrix} dx \, dy + \tfrac{1}{2}\iint_{G\times G} \tilde{r}^T(x,y) \begin{bmatrix} w(x) \\ v(y) \end{bmatrix} dx \, dy \, ,$$

--- (2.4)

where

$$\tilde{K}(x,y) := \begin{bmatrix} \dfrac{K_1(x)}{|G|} & K_2(x,y) \\ K_2(y,x) & \dfrac{K_1(y)}{|G|} \end{bmatrix}, \quad |G| := \int_G dx \; ; \quad \tilde{r}(x,y) := \dfrac{1}{|G|}\begin{bmatrix} r(x) \\ r(y) \end{bmatrix} \; .$$

--- (2.5)

We note that $\tilde{K}$ is a symmetric matrix, since $K_1$ is a symmetric matrix and $K_2(y,x) = K_2^T(x,y)$.

Then it can be verified that $E(w) = \tilde{E}(w,w)$.

We are interested in conditions under which the purely quadratic part of $E$ is a positive definite operator, which means that $E(w) - \int_G r_0^T(x) w(x) dx > 0$ where $w$ is any nonzero element of $L^2(G \mapsto \mathbb{R}^n)$.

The reason is that the positive definiteness of the purely quadratic part of $E$ is a sufficient condition for a critical point of $E$ to be a point of minimum.

We shall denote by $E_q$ the purely quadratic part of $E$, i.e.

$$E_q(w) := \tfrac{1}{2} \int_G w^T(x) K_1(x) w(x) dx + \tfrac{1}{2} \iint_{G \times G} w^T(x) K_2(x,y) w(y) dx dy \qquad \text{--- (2.6)}$$

and similarly for $\tilde{E}$, i.e.

$$\tilde{E}_q(w,v) := \tfrac{1}{4} \iint_{G \times G} [w^T(x) \ v^T(y)] \tilde{K}(x,y) \begin{bmatrix} w(x) \\ v(y) \end{bmatrix} dx dy \qquad \text{--- (2.7)}$$

Definition 2.1. We shall say that the pair of matrix-valued kernels ($K_1$, $K_2$) is a <u>pair that generates a positive-definite integral form</u> if $K_1(x)$ is invertible for all $x$ in $\overline{G}$ and, for every nonzero $w$ in $\left(L^2(G \mapsto \mathbb{R}^n)\right)$, we have

$$\tfrac{1}{2} \int_G w^T(x) K_1(x) w(x) dx + \tfrac{1}{2} \iint_{G \times G} w^T(x) K_2(x,y) w(y) dx dy > 0 \quad . \ ///$$

We can characterize this property of a pair of kernels in a sharp way that involves only one matrix-valued function.

When $K_1(x)$ is a positive-definite matrix for all $x$ in $\overline{G}$, we utilize the positive-definite square root of $K_1(x)$, which we denote by $K_1^{\frac{1}{2}}(x)$, and we set $v(x) := K_1^{\frac{1}{2}}(x) w(x)$. Then

$$w^T(x) K_1(x) w(x) + \int_G w^T(x) K_2(x,y) w(y) dx dy =$$
$$= \|v\|^2_{L^2(G)} + \int_G v^T(x) K_1^{-\frac{1}{2}}(x) K_2(x,y) K_1^{-\frac{1}{2}}(y) v(y) dx dy \ . \qquad \text{--- (2.8)}$$

Thus the condition that the pair ($K_1$, $K_2$) generates a positive-definite integral form is equivalent to the condition that, for all nonzero $v$ in $L^2(G)$, we have

$$\|v\|^2_{L^2(G)} + \int_G v^T(x) K_1^{-\frac{1}{2}}(x) K_2(x,y) K_1^{-\frac{1}{2}}(y) v(y) dx dy > 0 \quad . \qquad \text{--- (2.9)}$$

This is one type of <u>coercivity condition</u>, it means that the bilinear integral form generated by the kernel $K_1^{-\frac{1}{2}}(x) K_2(x, y) K_1^{-\frac{1}{2}}(y)$ satisfies

$$\int_G v^T(x) K_1^{-\frac{1}{2}}(x) K_2(x, y) K_1^{-\frac{1}{2}}(y) v(y) \, dx \, dy \geq -\|v\|_{L^2(G)}^2 \quad \text{for all } v, \text{ with equality only when } v \equiv 0.$$

This coercivity condition will be satisfied if $K_1^{-\frac{1}{2}}(x) K_2(x, y) K_1^{-\frac{1}{2}}(y)$ has the form of an expansion akin to a Mercer expansion of positive definite kernels, namely

$$K_1^{-\frac{1}{2}}(x) K_2(x, y) K_1^{-\frac{1}{2}}(y) = \sum_{k=1}^{\infty} \lambda_k w_k(x) w_k^T(y)$$

where $\lambda_k > -1 \ \forall k \in \mathbb{N}$, $\sum_{k=1}^{\infty} |\lambda_k| < \infty$, and $\{w_k\}_{k=1}^{\infty}$ is an orthonormal basis of $L^2(G \mapsto \mathbb{R}^n)$.

We shall prove

<u>Theorem 2.1.</u> When $K_1(x)$ is positive definite for all $x$ in $G$ and the pair $(K_1, K_2)$ is a pair that generates a positive-definite integral form, then:

(i) The Fredholm integral equation

$$K_1(x) w(x) + \int_G K_2(x, y) w(y) \, dy + r_0(x) = 0 \qquad \text{--- (2.10)}$$

has a unique solution $w^*$;

(ii) The unique solution $w^*$ minimizes $E(w)$ over all $w$ in $\left(L^2(G \mapsto \mathbb{R}^n)\right)$.

<u>Proof.</u> For proving assertion (i), we observe that Eq. (2.10) can be written as

$$w(x) = -\int_G K_1^{-1}(x) K_2(x, y) w(y) \, dy - K_1^{-1}(x) r_0(x) \qquad \text{--- (2.11)}$$

which is the standard form of a second-kind Fredholm integral equation, and Fredholm's alternative implies that either (2.11) has a unique solution, or the corresponding homogeneous equation

$$w(x) = -\int_G K_1^{-1}(x) K_2(x, y) w(y) dy$$

has a nonzero solution. If $w_1$ is a nonzero solution of the homogeneous equation, we have

$$K_1(x) w_1(x) + \int_G K_2(x, y) w_1(y) dy = 0$$

thus, by multiplying from the left by $w_1^T(x)$ and integrating over $G$, we obtain

$$\int_G w_1^T(x) K_1(x) w_1(x) dx + \iint_{G \times G} w_1^T(x) K_2(x, y) w_1(y) dy = 0 \ ,$$

which contradicts the positivity of $E_q$.

We now prove assertion (ii). Let $w^*$ be the unique solution of (2.10). Then, from

$$K_1(x) w^*(x) + \int_G K_2(x, y) w^*(y) dy = -r_0(x)$$

it follows that

$$E_q(w^*) = \tfrac{1}{2}\int_G w^{*T}(x) K_1(x) w^*(x) dx + \tfrac{1}{2}\iint_{G \times G} w^{*T}(x) K_2(x, y) w^*(y) dy = -\tfrac{1}{2}\int_G r_0^T(x) w^*(x) dx \ .$$

Consequently, we have $E_q(w^*) = -E(w^*)$ , because

$$E(w^*) = E_q(w^*) + \int_G r_0^T(x) w^*(x) dx = -\tfrac{1}{2}\int_G r_0^T(x) w^*(x) dx + \int_G r_0^T(x) w^*(x) dx =$$

$$= \tfrac{1}{2}\int_G r_0^T(x) w^*(x) dx = -E_q(w^*) \ .$$

--- (2.12)

We calculate $E_q(w - w^*)$:

$$E_q(w - w^*) = \tfrac{1}{2} \int_G (w(x) - w^*(x))^T K_1(x)(w(x) - w^*(x))\,dx +$$
$$+ \tfrac{1}{2} \iint_{G \times G} (w(x) - w^*(x))^T K_2(x,y)(w(y) - w^*(y))\,dx\,dy =$$
$$= E_q(w) + E_q(w^*) - \int_G w^T(x) K_1(x) w^*(x)\,dx - \iint_{G \times G} w^T(x) K_2(x,y) w^*(y)\,dx\,dy =$$
$$= E_q(w) + E_q(w^*) + \int_G r_0^T(x) w(x)\,dx = E(w) + E_q(w^*) = E(w) - E(w^*)\ .$$

For $w \neq w^*$, we have, by the positive definiteness condition on the pair $(K_1, K_2)$,

$$E_q(w - w^*) > 0$$

thus

$$E(w^*) < E(w)$$

which shows that $w^*$ is the unique minimizer of $E$. ///

We note that eq. (2.12) has independent significance: it gives an expression for the minimum value of $E$, namely

$$E(w^*) = \tfrac{1}{2} \int_G r_0^T(x) w^*(x)\,dx\ .$$

Also, we have:

<u>Corollary 2.1.</u> Under the same conditions for the kernel $K_1$, but changing the positive definiteness condition on the pair $(K_1, K_2)$ to a condition of positive semi-definiteness, namely

$$\tfrac{1}{2} \int_G w^T(x) K_1(x) w(x)\,dx + \tfrac{1}{2} \iint_{G \times G} w^T(x) K_2(x,y) w(y)\,dx\,dy \geq 0$$

for all $w$ in $\left(L^2(G \mapsto \mathbb{R}^n)\right)$, we obtain the conclusion that the Fredholm linear integral equation (2.10) is a necessary condition that every minimizer of $E(w)$ must satisfy. ///

3. Quadratic control of linear Fredholm integral equations.

We consider a controlled Fredholm integral equation

$$\varphi(x) = \varphi_0(x) + \int_G \{A(x,y)\varphi(y) + B(x,y)u(y)\}dy.$$

--- (3.1)

We shall study the optimal control problem of minimizing a quadratic cost functional of the form

$$J := \int_G \{\tfrac{1}{2}\varphi^T(x)P(x)\varphi(x) + \varphi^T(x)Q(x)u(x) + \tfrac{1}{2}u^T(x)R(x)u(x)\}dx.$$

--- (3.2)

The solution of (2.1), under appropriate assumptions, can be represented in terms of a resolvent kernel $K(x,y)$ as follows:

$$\varphi(x) = \varphi_0(x) + \int_G K(x,y)\varphi_0(y)dy + \int_G B(x,y)u(y)dy + \int_G \int_G K(x,z)B(z,y)u(y)dy\,dz.$$

--- (3.3)

The resolvent kernel $K$ satisfies

$$K(x,y) = A(x,y) + \int_G A(x,z)K(z,y)dz.$$

--- (3.4)

Eq. (3.3) allows us to represent the solution of (3.1) in the form

$$\varphi(x) = \varphi_1(x) + \int_G B_1(x,y)u(y)dy,$$

--- (3.5)

where

$$\varphi_1(x) := \varphi_0(x) + \int_G K(x,y)\varphi_0(y)dy, \quad B_1(x,y) := B(x,y) + \int_G K(x,z)B(z,y)dz.$$

--- (3.6)

The representation (3.5) transforms the cost functional $J$ into an integral functional that is quadratic in the control:

$$J = \tfrac{1}{2}\int_G \varphi_1^T(x)P(x)\varphi_1(x)\,dx + \int_G \int_G \varphi_1^T(y)P(x)B_1(y,x)u(x)\,dx\,dy +$$

$$+ \tfrac{1}{2}\int_G \int_G \int_G u^T(y)B_1^T(z,y)P(z)B_1(z,x)u(x)\,dx\,dy\,dz + \int_G \varphi_1^T(x)Q(x)u(x)\,dx +$$

$$+ \int_G \int_G u^T(y)B_1^T(x,y)Q(x)u(x)\,dx\,dy + \tfrac{1}{2}\int_G u^T(x)R(x)u(x)\,dx$$

--- (3.7)

from which the variation $\delta J$ of the cost functional $J$, under a variation $\delta u$ of the control, is calculated as

$$\delta J = \int_G \int_G \varphi_1^T(y)P(x)B_1(y,x)\,\delta u(x)\,dx\,dy +$$

$$+ \int_G \int_G \int_G u^T(y)B_1^T(z,y)P(z)B_1(z,x)\,\delta u(x)\,dx\,dy\,dz + \int_G \varphi_1^T(x)Q(x)\,\delta u(x)\,dx +$$

$$+ \int_G \int_G u^T(y)[B_1^T(x,y)Q(x) + Q^T(y)B_1(y,x)]\delta u(x)\,dx\,dy + \int_G u^T(x)R(x)\,\delta u(x)\,dx \;.$$

--- (3.8)

The vanishing of $\delta J$ yields the following integral equation for a stationary point of the functional $J$:

$$u^*(x) = -R^{-1}(x)\Bigg[ Q^T(x)\varphi_1(x) + \int_G B_1^T(z,x)P(x)\varphi_1(z)\,dz +$$

$$+ \int_G \bigg\{ [B_1^T(x,y)Q(x) + Q^T(y)B_1(y,x)] + \int_G B_1^T(z,x)P(z)B_1(z,y)\,dz \bigg\} u^*(y)\,dy \Bigg] \;.$$

--- (3.9)

This is a Fredholm integral equation of the second kind that is linear in the unknown function $u^*$. Consequently, standard theorems about existence and uniqueness of solutions of linear Fredholm integral equations apply.

Sufficient conditions, both for the existence of solutions of (3.8), and for those solutions to actually provide a minimum of the original cost functional, can be found by reducing the problem to minimizing a quadratic functional, and then invoking the results of section 2. Substitution of (3.5) into (3.2) shows that $J$ has the form

$$J = \int_G \{\tfrac{1}{2} u^T(x) K_1(x) u(x) + r^T(x) u(x)\} dx + \iint_{G \times G} \tfrac{1}{2} u^T(x_1) K_2(x_1, x_2) u(x_2) dx_2 \, dx_1 +$$
$$+ \int_G \{\tfrac{1}{2} \varphi_1^T(x) P(x) \varphi_1(x) dx$$

--- (3.10)

where

$$K_1(x) := R(x);$$
$$K_2(x_1, x_2) := \int_G B_1^T(y, x_1) R(y) B_1(y, x_2) dy + B_1^T(x_2, x_1) Q(x_2) + Q^T(x_1) B_1(x_1, x_2);$$
$$r^T(x) := \int_G \varphi_1^T(y) Q(x) B_1(y, x) dy + \varphi_1^T(x) Q(x) \ .$$

--- (3.11)

Consequently, sufficient conditions for existence and uniqueness of the solution of (3.9), and sufficient conditions for that solution to give a minimum of (3.2), are that (i) $R(x)$ be positive definite for all *x* in *G*, and (ii) that the pair $(R, K_2)$ generate a positive definite quadratic form on $L^2((G \times G) \to \mathbb{R}^m)$

.

## 4. Control of certain nonlinear Fredholm integral equations.

Now we consider a controlled Fredholm integral equation that is nonlinear in the state but linear in the control:

$$\varphi(x) = \varphi_0(x) + \int_G \{f(x,y,\varphi(y)) + F(x,y,\varphi(y))u(y)\}dy \qquad \text{--- `(4.1)}$$

with corresponding cost functional that is nonlinear in the state and quadratic in the control:

$$J := \int_G \{g_0(x,\varphi(x)) + g_1^T(x,\varphi(x))u(x) + \tfrac{1}{2}u^T(x)G(x,\varphi(x))u(x)\}dx. \qquad \text{--- (4.2)}$$

In this case, we proceed via an extremum principle. The Hamiltonian is a functional of the co-state $\psi$, and it is given by

$$H(x,\varphi,u,\psi(\cdot)) = g_0(x,\varphi) + g_1^T(x,\varphi)u + \tfrac{1}{2}u^T G(x,\varphi)u + \int_G \psi(y)[f(y,x,\varphi) + F(y,x,\varphi)u]dy .$$

$$\text{--- (4.3)}$$

The co-state $\psi$ satisfies $\psi(x) = \nabla_\varphi H(x,\varphi,u,\psi(\cdot))$, that is,

$$\psi(x) = \nabla_\varphi g_0(x,\varphi) + \nabla_\varphi g_1^T(x,\varphi)u + \tfrac{1}{2}u^T \nabla_\varphi G(x,\varphi)u + \int_G \psi(y)[\nabla_\varphi f(y,x,\varphi) + \nabla_\varphi F(y,x,\varphi)u]dy .$$

$$\text{--- (4.4)}$$

An optimal control $u^*$, together with the corresponding state function $\varphi^*$ and the corresponding co-state $\psi^*$, satisfy $\nabla_u H(x,\varphi^*(x),u^*(x),\psi^*(\cdot)) = 0$, that is,

$$g_1^T(x,\varphi^*(x)) + u^{*T}G(x,\varphi^*(x)) + \int_G \psi^*(y)F(y,x,\varphi^*(x))dy = 0$$

$$\text{--- (4.5)}$$

from which we find

$$u^*(x) = -G^{-1}(x,\varphi^*(x))\left[ g_1(x,\varphi^*(x)) + \int_G F^T(y,x,\varphi^*(x))\psi^{*T}(y)dy \right]. \qquad \text{--- (4.6)}$$

Substitution of (4.6) into (4.4) gives the following double Fredholm integral equation for an optimal co-state:

$$\psi^*(x) = \nabla_\varphi g_0(x, \varphi^*(x)) - \nabla_\varphi g_1^T(x, \varphi^*(x)) G^{-1}(x, \varphi^*(x)) \cdot$$

$$\cdot \left[ g_1(x, \varphi^*(x)) + \int_G F^T(y, x, \varphi^*(x)) \psi^{*T}(y) dy \right] +$$

$$+ \tfrac{1}{2} \iint_{G \times G} \psi^*(y) F(y, x, \varphi^*(x)) G^{-1}(x, \varphi^*(x)) \nabla_\varphi G(x, \varphi) G^{-1}(x, \varphi^*(x)) F^T(z, x, \varphi^*(x)) \psi^{*T}(z) dz\, dy +$$

$$+ \tfrac{1}{2} g_1^T(x, \varphi^*(x)) G^{-1}(x, \varphi^*(x)) \nabla_\varphi G(x, \varphi) G^{-1}(x, \varphi^*(x)) g_1(x, \varphi^*(x)) +$$

$$+ g_1^T(x, \varphi^*(x)) G^{-1}(x, \varphi^*(x)) \nabla_\varphi G(x, \varphi) \int_G F^T(y, x, \varphi^*(x)) \psi^{*T}(y) dy +$$

$$+ \int_G \psi^*(y) [\nabla_\varphi f(y, x, \varphi^*(x)) - \nabla_\varphi F(y, x, \varphi^*(x)) G^{-1}(x, \varphi^*(x)) g_1(x, \varphi^*(x))] dy -$$

$$- \iint_{G \times G} \psi^*(y) \nabla_\varphi F(y, x, \varphi) G^{-1}(x, \varphi^*(x)) F^T(z, x, \varphi^*(x)) \psi^{*T}(z) dz\, dy .$$

--- (4.7)

Also, substitution of (4.6) into (4.1) yields

$$\varphi^*(x) = \varphi_0(x) + \int_G \{ f(x, y, \varphi^*(y)) - G^{-1}(y, \varphi^*(y)) g_1(y, \varphi^*(y)) \} dy -$$

$$- \iint_{G \times G} F(x, y, \varphi^*(y)) G^{-1}(y, \varphi^*(y)) F^T(z, y, \varphi^*(y)) \psi^{*T}(z) dz\, dy .$$

--- (4.8)

The system of (4.7) and (4.8) consists of double Fredholm integral equations. This system is of the second degree in the co-state $\psi$, and in that respect it resembles the well-known Riccati differential system that arises in linear-quadratic optimal control for systems governed by differential equations.`        `

## 5. Quadratic control for linear Volterra controlled systems.

We consider a Volterra controlled system of the form

$$y(t) = y_0(t) + \int_0^t \{A(t,s)y(s) + B(t,s)u(s)\}\,ds \quad . \qquad \text{--- (5.1)}$$

The state $y$ is $n$-dimensional, the control $u$ is $m$-dimensional. $A$ and $B$ are matrices of dimensions compatible with the dimensions of the state and the control.

The control objective is to minimize a quadratic functional

$$J := \int_0^T \{\tfrac{1}{2} y^T(t)P(t)y(t) + y^T(t)Q(t)u(t) + \tfrac{1}{2} u^T(t)R(t)u(t)\}\,dt \quad . \qquad \text{--- (5.2)}$$

The Hamiltonian $H$ for this case is

$$H(t, y, \psi(\cdot), u) := \tfrac{1}{2} y^T P(t) y + y^T Q(t) u + \tfrac{1}{2} u^T R(t) u + \int_t^T \psi(s)\{A(s,t)y + B(s,t)u\}\,ds \quad .$$

$$\text{--- (5.3)}$$

The co-state $\psi$ satisfies $\psi(t) = \nabla_y H(t, y, \psi(\cdot), u)$, thus

$$\psi(t) = y^T P(t) + u^T Q^T(t) + \int_t^T \psi(s)A(s,t)\,ds \quad . \qquad \text{--- (5.4)}$$

If we denote by $u_*$ a minimizer of the Hamiltonian, for every $t, y, \psi(\cdot)$, we find (by differentiating the Hamiltonian with respect to $u$)

$$y^T(t)Q(t) + u_*^T(t)R(t) + \int_t^T \psi(s)B(s,t)\,ds = 0 \quad . \qquad \text{--- (5.5)}$$

By substituting the solution, for $u_*^T(t)$ into (5.4), we find the following equation for the co-state $\psi$:

$$\psi(t) = y^T(t)[P(t) - Q(t)R^{-1}(t)Q^T(t)] + \int_t^T \psi(s)[A(s,t) - B(s,t)R^{-1}(t)Q^T(t)]\,ds \quad .$$

$$\text{--- (5.6)}$$

Also, by substituting $u_*(t)$ into (5.1) and taking transposition of the resulting integral equation, we obtain

$$y^T(t) = y_0^T(t) + \int_0^t y^T(s)[A^T(t,s) - Q(s)R^{-1}(s)B^T(t,s)]ds - \int_{s=0}^T \int_{\sigma=s}^T \psi(\sigma)B(\sigma,t)R^{-1}(s)B^T(t,s)\,d\sigma\,ds$$

--- (5.7)

which, by a change of the order of integration in the double integral, becomes

$$y^T(t) = y_0^T(t) +$$
$$+ \int_0^t y^T(s)[A^T(t,s) - Q(s)R^{-1}(s)B^T(t,s)]ds - \int_{\sigma=0}^T \int_{s=0}^{\min(t,\sigma)} \psi(\sigma)B(\sigma,t)R^{-1}(s)B^T(t,s)\,d\sigma\,ds \,.$$

--- (5.8)

We set

$$K_1(t,\sigma) := \int_{s=0}^{\min(t,\sigma)} B(\sigma,t)R^{-1}(s)B^T(t,s)\,ds \quad;$$
$$C(t,s) := A^T(t,s) - Q(s)R^{-1}(s)B^T(t,s) \,.$$

--- (5.9)

Then the equation for the state becomes

$$y^T(t) = y_0^T(t) + \int_0^t y^T(s)C(t,s)\,ds - \int_0^T \psi(s)K_1(t,s)\,ds \,.$$

--- (5.10)

Further, let $S(t,\sigma)$ be the resolvent kernel associated with the kernel $C(t,s)$ in (5.10), and define

$$S_1(t,\sigma) := \int_{s=\sigma}^t K_1(s,\sigma)S(t,s)\,ds \,.$$

--- (5.11)

Then the state admits this representation:

$$y^T(t) = y_0^T(t) + \int_0^t y_0^T(s)S_1(t,s)\,ds - \int_0^T \psi(s)K_1(t,s)\,ds - \int_0^t \psi(s)S_1(t,s)\,ds \,.$$

--- (5.12)

Substitution of (5.12) into equation (5.6) for the co-state $\psi$ gives the final form of a second-kind Fredholm integral equation for $\psi$:

$$\psi(t) = \left[ y_0^T(t) + \int_0^t y_0^T(s)\,ds - \int_0^T \psi(s)K_1(t,s)\,ds - \int_0^t \psi(s)S_1(t,s)\,ds \right] \cdot$$
$$\cdot \left[ P(t) - Q(t)R^{-1}(t)Q^T(t) \right] + \int_t^T \psi(s)[A(s,t) - B(s,t)R^{-1}(t)Q^T(t)]ds \,.$$

--- (5.13)

Assuming the existence of a minimizing control, the determination of an optimal control depends on the solution of the Fredholm-Volterra integral equation (5.13). If a solution $\psi$ has been found, then tracing back through the calculations in this section, the corresponding optimal state trajectory is found from (5.10) or from (5.12), and the optimal control yielding that state trajectory is found from (5.5) which is a linear algebraic equation for that optimal control.